\documentclass[11pt]{amsart}
\usepackage{latexsym}
\usepackage{amsfonts}

\usepackage{amsmath,amsthm,amssymb}

\title[Generic properties]{Generic properties of Whitehead's algorithm
and isomorphism rigidity of random one-relator groups}

\author[I.~Kapovich]{Ilya Kapovich}

\address{\tt Department of Mathematics, University of Illinois at
  Urbana-Champaign, 1409 West Green Street, Urbana, IL 61801, USA
\newline http://www.math.uiuc.edu/\~{}kapovich/}
\email{\tt kapovich@math.uiuc.edu}

\author[P.~Schupp]{Paul Schupp}

\address{\tt Department of Mathematics, University of Illinois at
  Urbana-Champaign, 1409 West Green Street, Urbana, IL 61801, USA
\newline http://www.math.uiuc.edu/People/schupp.html}
\email{schupp@math.uiuc.edu}

\author[V.~Shpilrain]{Vladimir Shpilrain} \address{\tt Department of Mathematics, The City College of New York, New York, NY 10031, USA
\newline http://www.sci.ccny.cuny.edu/\~{}shpil}
\email{\tt shpil@groups.sci.ccny.cuny.edu}

\newtheorem{theor}{Theorem}

\newtheorem{thm}{Theorem}[section]
\newtheorem{lem}[thm]{Lemma}
\newtheorem{cor}[thm]{Corollary}

\newtheorem{prop}[thm]{Proposition}
\theoremstyle{definition}
\newtheorem{defn}[thm]{Definition}

\newtheorem{conv}[thm]{Convention}
\newtheorem{rem}[thm]{Remark}
\newtheorem{exmp}[thm]{Example}

\begin{document}

\begin{abstract}
  We prove that  Whitehead's algorithm for solving
  the automorphism problem in a fixed free group $F_k$ 
  has strongly linear time generic-case complexity. This is done by
  showing that the ``hard'' part of the algorithm terminates in
  linear time  on an exponentially   generic set of input pairs.
  We then apply these results to one-relator groups.
  We   obtain a  Mostow-type isomorphism rigidity result
  for random one-relator groups:  If two such groups are
  isomorphic then their Cayley graphs on the \emph{given generating sets} are
  isometric.  Although no nontrivial examples were previously
  known, we  prove that one-relator groups are generically
  \emph{complete} groups, that is, they have trivial center and trivial outer
  automorphism group.    We also prove that the
  stabilizers of generic elements of $F_k$ in $Aut(F_k)$ are cyclic
  groups generated by inner automorphisms and that $Aut(F_k)$-orbits
  are uniformly small in the sense of their growth entropy. We further
  prove  that the number $I_k(n)$ of
  \emph{isomorphism types} of $k$-generator one-relator groups with
  defining relators of length $n$ satisfies
\[
\frac{c_1}{n} (2k-1)^n \le I_k(n)\le \frac{c_2}{n} (2k-1)^n,
\]
where $c_1=c_1(k)>0, c_2=c_2(k)>0$ are some constants independent
of $n$.  Thus  $I_k(n)$ grows in essentially the same manner as the
  number of cyclic words of length $n$.
\end{abstract}


\maketitle

\tableofcontents

\section{Introduction}\label{intro}

The famous Mostow Rigidity Theorem~\cite{Mos73} says that if $M_1$ and 
$M_2$ are complete connected hyperbolic manifolds of finite volume
and dimension $n\ge 3$ then their fundamental groups are
isomorphic if and only if the manifolds themselves are isometric.
For a finitely generated group $G$ with a finite generating set
$A$ the naturally associated geometric object is the Cayley graph
$\Gamma(G,A)$. Thus one might say that a class of groups 
equipped with specified finite generating sets has the  \emph{isomorphism rigidity property} if
whenever two groups from this class are isomorphic then their
Cayley graphs on  the \emph{given} generating sets
are isometric. Phenomena of this  type were known for various
classes of Coxeter and Artin groups (e.g \cite{Ros,PrSp,Bahls,BMMN,MW}).
 In the present paper we obtain the first result of
this kind  for a class of groups given in terms of "general" finite
presentations. We prove if  two "random" one-relator groups
$G_u=\langle a_1,\dots, a_k |u=1\rangle$ and $G_v=\langle
a_1,\dots, a_k |v=1\rangle$ are isomorphic  then their  Cayley
graphs $\Gamma(G_u, \{a_1,\dots, a_k\})$ and $\Gamma(G_v,
\{a_1,\dots, a_k\})$ are isometric.  Indeed, their
Cayley graphs are  isomorphic as labeled graphs 
by a  graph isomorphism which is only allowed to
permute the label set $\{a_1,\dots,
a_k\}^{\pm 1}$. This provides a conceptually new source of
group-theoretic rigidity given by "random" or "generic" groups. 
Such  rigidity arises not from 
structural restrictions, such as the structure of flats or
of finite subgroups,  but rather  from the rigidity  of
"randomness''  itself. 

The theorems in  this paper are based on combining  very
different probabilistic and algebraic techniques: the  generic-case
analysis of Whitehead's algorithm in this paper and earlier  work of
Kapovich-Schupp~\cite{KS} on the Nielsen Uniqueness property for
generic groups that utilized the Arzhantseva-Ol'shansakii~\cite{AO} graph minimization and genericity techniques.  Our  goal is to obtain
new algebraic and geometric applications and the 
probabilistic tool used in this paper,  Large Deviation Theory
applied to finite state Markov chains, is  quite  basic from the point
of view of probability theory. Nevertheless, combining it  with
algebraic and algorithmic considerations  as well as with earlier
probabilistic results on  Nielsen Uniqueness produces
surprisingly powerful results.

We adopt the following convention throughout this paper.

\begin{conv}\label{conv:main}
Let   $F_k=F(a_1,\dots, a_k)$ be the free group of rank  $k\ge 2$. 
The \emph{group alphabet} is  $\Sigma:=\{a_1,\dots,a_k\}^{\pm 1}$.
 A word $w\in \Sigma^*$ is \emph{reduced} if $w$ does not contain 
any subwords of the   form $a_ia_i^{-1}$ or $a_i^{-1}a_i$. 
The \emph{length}, $|w|$, of a word $w$ is the number of letters in $w$.
Since every element of $F_k$ can be represented by a unique
  reduced word, we can identify elements of $F_k$ with reduced words.
The length $|g|$ of an element element $g\in F_k$ 
is the length of the unique  reduced word  in $\Sigma^*$ which represents  $g$.

A word $w$ is \emph{cyclically  reduced} if all cyclic permutations of $w$ are reduced. 
We use $C$ to denote the set of all   cyclically reduced words in $F_k$.
Any   reduced word  $w$ can be uniquely decomposed as a concatenation
 $w   = vuv^{-1}$ where $u$ is a cyclically reduced.
  The   word $u$ is called the \emph{cyclically reduced form of $w$} and
  $||w||:=|u|$ is  the \emph{cyclic length} of $w$.

An element $w\in F_k$ is  \emph{minimal} if $|\phi(w)|\ge |w|$ 
 for all $\phi\in Aut(F_k)$. In other words, $w$ is a shortest
element  in its orbit $Aut(F_k)w$.
\end{conv}

Recall that the \emph{automorphism problem} (also called
the \emph{automorphic conjugacy problem} or the \emph{automorphic
  equivalence problem}) for a free group $F_k$  is the following decision problem:  Given two elements $u,v\in F_k$, is there an
automorphism $\phi\in Aut(F_k)$ such that $\phi(u)=v$?  If there is 
such an automorphism we say that $u$ and $v$ are \emph{automorphically equivalent}.
In a classic 1936 paper~\cite{Wh}
Whitehead provided an algorithm for solving this problem. 
We  need to give a  brief description of Whitehead's solution
 and more details are given in Section~\ref{w-moves} below.  Whitehead
introduced a particular finite set of generators  of $Aut(F_k)$, 
now called \emph{Whitehead automorphisms}. These
automorphisms are divided in two types. The Whitehead automorphisms of
the \emph{first kind} are ``relabeling automorphisms'' 
induced by permutations of the set $\{a_1,\dots, a_k\}^{\pm 1}$ and thus 
do not change the length of an element. The remaining Whitehead
automorphisms are of the \emph{second kind} and  can change the length of an
element. These automorphisms are precisely defined  in Definition~\ref{defn:moves} below.

\begin{prop}\label{wh}[Whitehead's  Theorem]~\cite{Wh}
\begin{enumerate}
\item (Length reduction) If $u \in F_k$ is cyclically reduced and not minimal  then there
  is a Whitehead automorphism $\tau$ such that $||\tau(u)||<||u||$.
\item (Length preservaton or ``peak reduction'') Let $u,v \in F_k$ be minimal 
 (and hence cyclically reduced)
  elements with $|u|=|v|=n>0$. Then $Aut(F_k)u=Aut(F_k)v$ if and only if
  there exists a finite sequence of Whitehead automorphisms
  $\tau_s,\dots, \tau_1$ such that $\tau_s\dots \tau_1(u)=v$ and such
  that for each $i=1,\dots, s$ we have
\[
||\tau_i\dots \tau_1(u)||=n.
\]

\end{enumerate}
\end{prop}

  This statement immediately gives Whitehead's  algorithm  for solving the automorphism
  problem for $F_k$.  First, by length reduction 
  there is a  algorithm 
  which, given any  element $w\in F_k$, finds a minimal element $w'\in
  Aut(F_k)w$. To start, cyclically reduce $w$. Then  repeatedly check
  if there is a Whitehead automorphism $\tau$ decreasing the
  cyclically reduced length of the current element and if so,
  apply such a $\tau$ and  cyclically reduce the result.
 This process terminates in at most $|w|$
steps with a minimal element and requires at worst quadratic time
in the length of $w$. Each step takes at most linear time since
the number of Whitehead automorphisms is fixed.
Thus given two elements of $F_k$ we can first replace them by
minimal $Aut(F_k)$-equivalent elements. By peak reduction, if these minimal elements
have different lengths  then there does not exist an automorphism
taking one of original elements to the other. This quadratic time
procedure is  the so-called ``easy part'' of Whitehead's algorithm.

   Now suppose that starting with  elements $u,v \in F_k$
the process above yields corresponding minimal elements $u',v'$
of the same length. Peak reduction  implies that if
these two minimal elements  are automorphically equivalent
then there is a chain of Whitehead automorphisms taking one element to
the other so that the \emph{cyclically reduced length is constant
  throughout the chain}. Since the number of elements of given length
is bounded by an exponential function, this provides an
algorithm  which is at worst exponential time
for deciding if two minimal elements of the same length are in the same
$Aut(F_k)$-orbit. This stage is called the ``hard part" of 
Whitehead's  algorithm.

 Taken together, these  two parts 
provide a complete solution for the automorphism problem for  $F_k$
and requires \emph{at most exponential time} in terms of the
maximum of the lengths of the input words. Note that Whitehead's
algorithm actually solves the \emph{Search Automorphism Problem}
as well. If $u,v$ are in the same $Aut(F_k)$-orbit, the
algorithm produces  an explicit automorphism taking $u$ to $v$.

   Whether or not Whitehead's algorithm actually requires
exponential time is currently an  active research question.
 The only well understood case is $k=2$ where Myasnikov and
Shpilrain~\cite{MS} proved that an improved version of Whitehead's
algorithm takes at most polynomial time.  Substantial further
progress for $k=2$ has been made by Bilal Khan~\cite{Khan}. Very interesting
partial results regarding the complexity of Whitehead's algorithm for
$k>2$ have recently been  obtained by Donghi Lee~\cite{Lee}.

Experimental evidence (for example~\cite{BB,HMM,MM}) strongly
indicates that even for $k>2$ Whitehead's algorithm usually runs
very quickly. In the present paper
we provide a theoretical explanation of this phenomenon and prove that
that for an ``exponentially generic'' set of inputs the ``easy'' first stage of
the Whitehead algorithm terminates immediately and the ``hard'' second
part terminates in linear time.

The study of genericity, or ``typical behavior'', in group theory was
initiated by Gromov~\cite{Grom,Grom1}, Ol'shanskii~\cite{Ol92} and
Champetier~\cite{Ch94}. The importance of these ideas is becoming
increasingly clear and manifestations of genericity in many different
group-theoretic contexts are the subject of active investigation
\cite{A1,A2,A3,AO,Ch94,Ch95,Ch00,Che96,Che98,Z,KS,Gh,Grom2,KMSS,KMSS1,Oliv}.
 Intuitively, a subset $Q$ of $S\subseteq F_k$ is generic in $S$ if a
``randomly'' chosen long element of $S$ belongs to $Q$ with
probability tending to $1$, or that $Q$ has ``measure 1'' in $S$.  
The precise definitions of genericity used  in~\cite{KMSS,KMSS1} are given in 
Definition~\ref{defn:generic} below.

We need the following crucial definition.
\begin{defn}
A cyclically reduced element $w \in F_k$ is  \emph{strictly minimal} if
 the cyclically reduced length  $||\tau(w)||$ is strictly greater than  $|w|$
 for every  non-inner Whitehead automorphism $\tau$ of the second kind.
 We use  $SM$ to denote the set of all strictly minimal elements of $F_k$.
  Also, $SM'$ denotes the set of all $w\in F_k$ such
that the cyclically reduced form of $w$ belongs to $SM$.
\end{defn}
  The description of Whitehead's algorithm given above 
shows  that every element of $SM$ is already
minimal in its $Aut(F_k)$-orbit.  Moreover, if $w\in SM$ then any
chain of Whitehead moves that preserves the cyclic length of
$w$ must consist entirely of conjugations and of Whitehead
automorphisms of the first kind, that is, relabeling automorphisms.
 Thus if $w\in SM$ and $w'\in F_k$ is
another minimal element with $|w|=|w'|$ then Whitehead's  algorithm,
applied to the pair $(w,w')$, terminates in time linear in $|w|$.
Moreover, for arbitrary $(w_1,w_2)\in F_k^2$ such that at least one of
$w_1, w_2$ is $Aut(F_k)$-equivalent to a strictly minimal element, then
Whitehead's  algorithm terminates in at most quadratic time on
$(w_1,w_2)$.

 We give here a short informal summary of our results
regarding Whitehead's algorithm and the properties of random
one-relator groups. Precise and detailed statements are given in Section~\ref{main}.
\begin{conv}
For $u \in F_k$ set $G_u = \langle a_1,...,a_k | u \rangle$.
\end{conv}
By saying that a certain
property holds for a generic element we mean that there is an
exponentially generic set such that every element of that  set has the
property. We prove that:

\begin{itemize}
\item[(a)] The cyclically reduced form of generic element of $F_k$
is strictly minimal and a generic cyclically reduced element is strictly minimal.
\item[(b)] The generic-case complexity of Whitehead's algorithm for $F_k$ 
   is  strongly linear-time.
\item[(c)] For any $u\in F_k$ the orbit $Aut(F_k)u$ is an
exponentially   negligible subset of $F_k$. Moreover, all such orbits are
  ``uniformly small'' in $F_k$. Namely, there is a number
  $\alpha<2k-1$ such that for any $u\in F_k$ the exponential growth
  rate of $Aut(F_k)u$ is $\le \alpha<2k-1$. (Note that the growth rate
  of $F_k$ is $2k-1$.).
\item[(d)] For a generic element $u\in F_k$ the stabilizer of $u$
  in   $Aut(F_k)$ is infinite cyclic and is generated by the inner
  automorphism corresponding to conjugation by $u$.
\item[(e)] For a generic $u\in F_k$ the one-relator group $G_u$ is
 a complete   group, that is, it has trivial center and trivial outer automorphism
  group.
\item[(f)] A generic one-relator group $G_u$ is torsion-free
  non-elementary   word-hyperbolic and it has either the Menger curve or the Sierpinski
  carpet as its  boundary. If $k=2$ the boundary is the Menger curve.
\item[(g)] If we \emph{fix} a generic one-relator group $G_u$ then
  there is   a quadratic-time algorithm (in terms of $|v|$) which decides if an
  \emph{arbitrary} one-relator group $G_v=\langle a_1,\dots, a_k
  |v\rangle$ is isomorphic to $G_u$.
\item[(h)] Two generic one-relator groups $G_u,G_v$ are isomorphic
if and only if $|u|=|v|$ and there is a relabeling automorphism
$\tau$ such that $\tau(u)$ is a cyclic permutation of $v$ or
$v^{-1}$.
 \item[(i)] The number $I_k(n)$ of \emph{isomorphism types} of
  one-relator groups on $k$ generators with defining relators of length $n$
  satisfies
\[
\frac{c_1}{n} (2k-1)^n \le I_k(n)\le \frac{c_2}{n} (2k-1)^n,
\]
where $c_1=c_1(k)>0, c_2=c_2(k)>0$ are some constants independent
of $n$.
\end{itemize}

The structure of Whitehead's algorithm for solving the
automorphism problem is  similar to that of Garside's
algorithm (and its various modifications) for solving the conjugacy
problem in braid groups. (See for example~\cite{Gar,BKL,FGM}.)
In both cases there has been
a great deal of experimental evidence that in practice the
algorithms almost always work much faster than the worst-case
exponential time estimate suggests. Statements (a) and (b) above
provide the first  proof explaining why this happens for
Whitehead's algorithm. It remains an interesting open problem to
find and prove similar statements for Garside's algorithm.

As discussed earlier, statement (h) above may be regarded as an analogue of
Mostow rigidity for random one-relator groups. 
Indeed, it says
that two generic one-relator groups $G_u$ and $G_v$ are isomorphic
if and only if their Cayley graphs corresponding to the
\emph{given} generating sets $\{a_1,\dots, a_k\}$ are isomorphic
as labelled graphs  where the graph isomorphism is only allowed to
permute the label set $\{a_1,\dots, a_k\}^{\pm
1}$. This means that the class of random one-relator groups has the isomorphism rigidity property.
 We will see that isomorphism rigidity is also responsible for us being
able to estimate  the number
of isomorphism types of one-relator groups in the statement (i)
above. In  subsequent work~\cite{KSn} Kapovich and Schupp combine
the results of this paper with  methods involving 
Kolmogorov complexity to prove that a random one-relator
presentation $G_u$ is ``essentially incompressible''. This means that $G_u$
does not admit any finite group presentation of total length much smaller than $|u|$.

\noindent{\bf Acknowledgements.}
We are grateful to Richard Sowers and Ofer Zeitouni for very
illuminating discussions regarding Large Deviation Theory.  We thank
Jean-Francois Lafont for raising the question of  counting the
number of isomorphism types of one-relator groups. We are also grateful to the referee for a number of comments that improved the paper.

\section{Generic sets and Generic Complexity}\label{Sect:generic}

We need to 
recall the definitions concerning  genericity used in~\cite{KMSS}. 
Note that the length condition on sets of pairs
which we consider here is slightly different from that used in~\cite{KMSS}.

We say that a  sequence $x_n\in \mathbb R$, $n\ge 1$ with $\lim_{n\to\infty}
  x_n=x\in \mathbb R$  \emph{converges  exponentially fast} if
  there are  $ 0<\sigma<1$ and $K>0$  such that for
  all $n\ge 1$ 

\[
|x-x_n|\le K \sigma^n.
\]

\begin{defn}\label{defn:generic}
  Let $S$ be a set of words in the  group alphabet $\Sigma$.
  Let  $\rho(n,S)$  denote 
  the number of words $w\in S$ with $|w|\le n$. Also, let
  $\gamma(n,S)$ denote the number of words $w\in S$ with $|w|=n$.

  We say that a subset $B\subseteq S$ is \emph{generic in $S$} if
\[
\lim_{n\to\infty} \frac{\rho(n,B)}{\rho(n,S)}=1.
\]
If, in addition, the convergence  is exponentially fast,
we say that $B$ is \emph{exponentially generic in $S$}.

The complement of an (exponentially) generic set in $S$ is said to be
\emph{(exponentially) negligible in $S$}.

Similarly, let $D\subset S\times S$ and let  $\rho(n,D)$ denote the number
of pairs $(u,v)\in D$ such that $|u|\le n$ and $|v|\le n$. Note that
$\rho(n, S\times S)= \rho(n,S)^2$.  We say that $D$ is \emph{generic
  in $S\times S$} if
\[
\lim_{n\to\infty} \frac{\rho(n,D)}{\rho(n,S\times S)}=1.
\]
Again, if  convergence is exponentially fast,
we say that $D$ is \emph{exponentially generic in $S\times S$}.
\end{defn}

  We can now apply this concept to decision problems. The following notion was introduced in \cite{KMSS}.

\begin{defn}[Generic-case complexity]\label{defn:gencom} 

  Let $S\subseteq \Sigma^*$ be an infinite set of words
  and let $D\subseteq S\times S$. (We regard the set $S\times S$ as the set of all inputs for a decision problem $D$, so that we are now working
relative to $S$).

Suppose that 
  $\Omega$ is a partial algorithm for deciding if an element $(u,v)\in
  S\times S$ belongs to $D$. Note that this means  that $\Omega$ is 
  \emph{correct}. That is, whenever
  $\Omega$ does produce a definite answer,  that answer is correct.
  Let $t(n)\ge 0$ be a non-decreasing function. We say that $\Omega$
  \emph{solves $D$ with strong  generic-case time complexity
    bounded by $t$ in $S\times S$} if there exists an exponentially $S\times
  S$-generic subset $A\subset S\times S$ such that for any $(u,v)\in
  A$ with $|u|\le n, |v|\le n$ the algorithm $\Omega$ terminates on
  the input $(u,v)$ in at most $t(n)$ steps.

  Let $S,D$ be as above and let $\mathcal B$ be a deterministic time
  complexity class such as linear time, quadratic time, polynomial time, etc. We say that $D$ is \emph{decidable with strong  $S$-generic
    case complexity in $\mathcal B$} if there exist a function $t(n)$
  satisfying the constraints of the complexity class $\mathcal B$ and a correct partial
  algorithm $\Omega$ that solves $D$ with strong  generic-case
  time complexity bounded by $t$ in $S\times S$.
\end{defn}

\section{Main results}\label{main}

We can now state   our main results regarding
Whitehead's algorithm in more technical detail.

\begin{theor}\label{thm:A}
  Let $F_k=F(a_1,\dots, a_k)$ where $k\ge 2$. Then
\begin{enumerate}
\item The set $SM\subseteq C$ is exponentially $C$-generic and the set
  $SM' \subseteq F_k$ is exponentially $F_k$-generic. Hence the set
  $SM\times SM\subseteq C\times C$ is exponentially $C\times
  C$-generic and the set $SM' \times SM' \subseteq F_k\times F_k$ is
  exponentially $F_k\times F_k$-generic.
\item There is a linear time (in $|w|$) algorithm which,
  given a freely reduced word $w$, decides whether or not $w\in SM$ and whether or not $w \in SM'$.
\item Every $w\in SM$ is minimal in its $Aut(F_k)$-orbit, that is for
  every  $\alpha\in Aut(F_k)$ we have $|w|\le |\alpha(w)|$.

  Moreover, if $w\in SM$ and $v$ is a cyclically reduced word with
  $|w|=|v|$ then $w$ and $v$ are in the same $Aut(F_k)$-orbit if and
  only if there exists a Whitehead automorphism $\tau$ of the first
  kind such that $\tau(w)$ is a cyclic permutation of $v$.
\item  Whitehead's  algorithm works in \emph{linear time} on pairs
  $(u,v)\in SM\times SM$ and so  has strongly
  linear time generic-case complexity on $C\times C$  Similarly,  Whitehead's
  algorithm works in linear time on pairs $(u,v)\in SM' \times SM'$
  and so has strongly linear time
  generic-case complexity on $F_k\times F_k$.
\item  Whitehead's  algorithm works in at most \emph{quadratic time}
  on \emph{all} pairs $(u,v)$ such that at least one of $u,v$ is in
  the same $Aut(F_k)$-orbit as an element of $SM$.
\end{enumerate}
\end{theor}

  The  theorem above says that for a ``random'' pair of cyclically
reduced words $(u,v)$ both $u$ and $v$ are strictly minimal. Hence the
``easy'' first part of Whitehead's algorithm  terminates in a
single step and the ``hard'' second  part  reduces to simply 
checking if one can get from  $u$ to  $v$ by applying a relabeling
automorphism and then a cyclic permutation.

Recall that for a subset $S\subseteq F_k$ the \emph{exponential growth
  rate} or \emph{growth entropy} of $S$ is
\[
H(S):=\limsup_{n\to\infty}\sqrt[n]{\rho(n,S)}.
\]
Then $H(F_k)=2k-1$ and $S\subseteq F_k$ is exponentially
$F_k$-negiligible if and only if $H(S)<2k-1$.

\begin{cor}\label{orbit}
  For any $w\in F_k$ the
  set $Aut(F_k)w$ is exponentially negligible in $F_k$ and the set
  $C\cap Aut(F_k)w$ is exponentially negligible in $C$. Moreover
\[
H(Aut(F_k)w)\le H(F-SM')<2k-1.
\]
\end{cor}
\begin{proof}
  We may assume that $w$ is minimal. Let $L$ be the set of elements of
  length $|w|$ in the orbit $Aut(F_k)w$. Now  $L$ is finite and any
  element in $Aut(F_k)w-L$ is not minimal  and hence  not strictly
  minimal. Therefore $T:=C\cap [Aut(F_k)w-L]\subseteq C-SM$. By
  part~(1) of Theorem~\ref{thm:A} the set $C-SM$ is exponentially
  $C$-negligible and therefore  so is the set $T$. 
  We have $C\cap Aut(F_k)w= T\cup (C\cap L)$ and therefore
  $C\cap Aut(F_k)w$ is $C$-negligible, as claimed.

  Let $u$ be an arbitrary element of $ Aut(F_k)w$.
  Since $u$ need not be   cyclically reduced let  $u_0$ 
  be the cyclically reduced form of   $u$.

  If $u_0 \notin SM$ then $u$ is contained in the set $F_k-SM'$
  which is exponentially $F_k$-negligible by part (1) of
  Theorem~\ref{thm:A}.  Now suppose that $u_0$ is strictly minimal. 
  Since $u_0$   is conjugate to $u$,  $u_0\in Aut(F_k)w$. Since $u_0$ is
  minimal, $|u_0|=|w|$ and $u_0\in L$. Thus $u$ is contained in the
  $F_k$-conjugacy class of an element of $L$. 
  It is not difficult to see that  any conjugacy class in $F_k$ 
  has  exponential growth rate  $\sqrt{2k-1}$ and is thus
  exponentially negligible. Therefore  the orbit
  $Aut(F_k)w$ is contained in the union of finitely many exponentially
  $F_k$-negligible sets and  is exponentially $F_k$-negligible, as required.

  Moreover, the set $F_k-SM'$
  contains the conjugacy class of $a_1$. Thus $H(F_k-SM')\ge
  \sqrt{2k-1}$.  The previous argument shows that $Aut(F_k)w$ is
  contained in the union of $F-SM'$ and of finitely many $F_k$-conjugacy
  classes $K_1, \dots, K_m$.  Hence

\begin{align*}
  &H(Aut(F_k)w)\le \\
  &\le\max\{ H(F_k-SM'), H(K_1),\dots, H(K_m)\}=H(F_k-SM')<2k-1,
\end{align*}
where the last inequality holds since $SM'$ is exponentially
$F_k$-generic and $F_k-SM'$ is exponentially $F_k$-negligible.
\end{proof}

Corollary~\ref{orbit} shows that automorphic orbits in $F_k$ are
``uniformly small'' in the sense of their growth rate. This can be
viewed as a generalization of the results of Borovik-Myasnikov-Shpilrain~\cite{BMS}
 and of Burillo-Ventura~\cite{BV} who established (with specific quantitative
growth estimates) that the set of primitive elements is exponentially
negligible in $F_k$.

As  mentioned before, the worst-case complexity of
Whitehead's algorithm is known to be polynomial time  for $k=2$.
The results of \cite{KMSS1} and Theorem~\ref{thm:A} imply that
the \emph{average-case} complexity  
(as opposed to  generic-case)  of Whitehead's algorithm is linear time for $k=2$.

A deep result of McCool~\cite{Mc} shows that for any $w\in F_k$ the
stabilizer of $w $ in $Aut(F_k)$ is finitely presentable.  Similar
arguments as those used in the proof of Theorem~\ref{thm:A} allow us to
conclude that $Aut(F_k)$-stabilizers of generic elements of $F_k$ are
very small.

\begin{defn}
  The set $TS$ (for ``Trivial Stabilizer'') is  the set of
  all words $w\in SM$  (necessarily cyclically reduced) such that $w$
  is not a proper power and such that for every nontrivial relabeling
  automorphism $\tau$ of $F_k$ the elements $w$ and $\tau(w)$ are not
  conjugate in $F_k$. Also,  $TS'$ denotes  the set of all elements of
  $F_k$ whose  cyclically reduced form is in $TS$.
\end{defn}

\begin{theor}\label{thm:B}
Let $k\ge 2$. Then:
\begin{enumerate}
\item The set $TS'$ is exponentially $F_k$-generic and the set $TS$ is
  exponentially $C$-generic.
\item There is a linear-time (in terms of $|w|$) algorithm which,
  given a freely reduced word $w$, decides if $w  \in TS'$ 
  or if  $w \in TS$.
\item For any nontrivial $w\in TS'$ the stabilizer $Aut(F_k)_w$ of $w$
  in $Aut(F_k)$ is the infinite cyclic group generated by the inner
  automorphism $ad(w)$ of $F_k$. Here 
  $ad(w):u\mapsto w u w^{-1}$ for  $u\in F_k$.
\item For every $w\in TS'$ the stabilizer $Out(F_k)_w$ of the conjugacy
  class of $w$ in $Out(F_k)$ is trivial.
\end{enumerate}
\end{theor}

  These  results  together with the  work of
Kapovich-Schupp~\cite{KS} on the isomorphism problem for one-relator
groups yield  strong conclusions about the properties of generic
one-relator groups. There are several different notions of genericity
in the context of finitely presented groups, namely genericity in the
sense of Arzhantseva-Ol'shanskii~\cite{AO} and in the sense of
Gromov~\cite{Grom,Ol92}. These two notions essentially coincide
in the case of one-relator groups.
Recall that a  group $G$ is \emph{complete} if all automorphisms of $G$ are inner (so that   $Out(G) = \{1\}$) and if $G$ also has trivial center so  that the adjoint map $ad:G\to Aut(G)$ is an isomorphism.

\begin{theor}\label{thm:C}
   There exists an exponentially
  $C$-generic set $Q_k$ of nontrivial cyclically reduced words in $F_k$ with
  the following properties:
\begin{enumerate}
\item There is an exponential time (in $|w|$) algorithm which, given a
  cyclically reduced word $w$, decides whether or not $w\in Q_k$.

\item Let $u\in Q_k$.  Then the one-relator group $G_u$ is a complete
  one-ended torsion-free word-hyperbolic group.

\item If $u\in Q_k$ then the hyperbolic boundary $\partial G_u$ is
  homeomorphic to either the Menger curve or the Sierpinski carpet. If
  $k=2$ then $\partial G_u$ is homeomorphic to the Menger curve.

\item Let $u,v\in Q_k$. Then the groups $G_u$ and $G_v$ are isomorphic
  if and only if there exists a relabeling automorphism $\tau$ of $F_k$
  such that $\tau(u)$ is a cyclic permutation of either $v$ or
  $v^{-1}$. In particular, $G_u\cong G_v$ implies $|u|=|v|$.

\item Let $u\in Q_k$ be a fixed element. Then there exists a quadratic
  time algorithm (in terms of $|v|$) which, given an \emph{arbitrary}
  $v\in F_k$, decides if the groups $G_u$ and $G_v$ are isomorphic.
\end{enumerate}
\end{theor}

It is worth noting that by a result of Champetier~\cite{Ch95},
obtained by completely different methods, generic (in the sense of
Gromov~\cite{Grom1,Ol92}) \emph{two-relator} groups are
word-hyperbolic with boundary homeomorphic to the Menger curve.

Prior to Theorem~\ref{thm:C} there were no known nontrivial examples
of complete one-relator groups and some experts in the field believed
that such groups might not exist. Our proof that such groups do exist
is obtained by an indirect probabilistic argument. The set $Q_k$ is
obtained as the intersection $Q_k=R_k\cap Z_k$ of two exponentially
$C$-generic sets,  $R_k$ and $Z_k$, and hence $Q_k$
is also exponentially generic.  In particular it is certainly non-empty.  
The genericity of the
sets $R_k$ and $ Z_k$ is established using two very different methods:
namely, the Arzhantseva-Ol'shanskii graph-minimization method 
 in~\cite{KS} and
Large Deviation Theory in the present paper. This
demonstrates the strength of the ``probabilistic argument''  for
producing groups with genuinely new and often unexpected features.

In the definitions of genericity both in the sense of
Gromov~\cite{Grom1,Ol92} and in the sense of Ol'shanskii~\cite{AO} one
counts group presentations as opposed to group isomorphism classes.
It is very natural to ask, for fixed numbers of generators and
defining relators, how many \emph{isomorphism types} there are of groups with
particular constraints on the lengths of the relators.
As a corollary of Theorem~\ref{thm:C}
it turns out that the number of \emph{isomorphism
  types} of one-relator groups with relators of length $n$ grows
in essentially the same way (taking into account the obvious
symmetries) as the number of one-relator presentations with relators
of length $n$.

\begin{cor}\label{classes}
  Let $k\ge 2$ be an integer. For $n\ge 1$ define $I_k(n)$ to be the number
  of \emph{isomorphism types} among the groups given by presentations
  $\langle a_1,\dots, a_k | u=1\rangle$ where $u$ varies of the set of
  all cyclically reduced words of length $n$. Then there exist
  constants $A=A(k)>0, B=B(k)>0$ such that for any $n\ge 1$
\[
\frac{B}{n} (2k-1)^n \le I_k(n)\le \frac{A}{n} (2k-1)^n.
\]
\end{cor}
\begin{proof}
  Let $Q_k$ be the exponentially generic set of cyclically reduced
  words given by Theorem~\ref{thm:C} and recall that 
  $C$ denotes the set of all cyclically reduced words.

  It follows from Lemma~\ref{crit} below that the number $\gamma(n,C)$
  of cyclically reduced words of length $n$ satisfies
\[
c_2 (2k-1)^n \ge \gamma(n,C)\ge c_1 (2k-1)^n
\]
for some constants $c_1,c_2>0$ independent of $n$.  

Since $Q_k$ is exponentially $C$-generic, Lemma~\ref{crit} below
implies that \[\lim_{n\to\infty} \frac{\gamma(n,Q_k)}{\gamma(n,C)}=1.\]
 Thus there is $n_0>1$ such
that for any $n\ge n_0$ we have
\[
\gamma(n,Q_k)\ge \frac{1}{2} \gamma(n,C)\ge \frac{c_1}{2} (2k-1)^n.
\]

Let $M$ be the number of all Whitehead automorphisms of the first kind
(that is, relabeling automorphisms).  Let $n\ge n_0$ and let $u\in
Q_k$ with $|u|=n$.  Part~4 of Theorem~\ref{thm:C} implies that the
number of $v\in Q_k$ with $G_v\cong G_u$ is $\le 2nM$. Here the factor
of $2n$ corresponds to the number of cyclic permutations of $u^{\pm
  1}$.

Therefore for $n\ge n_0$:
\[
I_k(n)\ge \frac{\gamma(n,Q_k)}{2Mn}\ge \frac{c_1}{4Mn} (2k-1)^n.
\]

The set $PP$ of cyclically reduced proper powers is exponentially
negligible in $C$ (see~\cite{AO}). Thus there exist $K>0$ and
$0<\sigma<1$ such that for any $n\ge 1$ we have
\[
\gamma(n,PP)\le K \sigma^n \gamma(n,C)\le Kc_2 \sigma^n (2k-1)^n.
\]

It is easy to see that if $u$ is cyclically reduced of length $n$ and
is not a proper power, then all $n$ cyclic permutations of $u$ are
distinct words. Clearly, if $v$ is a cyclic permutation of $u$ then
$G_u\cong G_v$.

Therefore

\[
I_k(n)\le \frac{\gamma(n,C-PP)}{n}+\gamma(n,PP)\le \frac{c_2}{n}
(2k-1)^n+\gamma(n,PP)\le \frac{2c_2}{n} (2k-1)^n,
\]
where the last inequality holds for all sufficiently large $n$.
\end{proof}

 Via an additional technical argument, Kapovich and Schupp~\cite{KSn} improve the  estimate for $I_k(n)$ and establish that
\[ \lim_{n\to\infty} \frac{n I_k(n)}{(2k-1)^n}=\frac{1}{k!2^{k+1}}. \]

\section{Whitehead automorphisms }\label{w-moves}

We follow Lyndon and Schupp, Chapter~I~\cite{LS} in 
 recalling  the basic definitions and results about
Whitehead automorphisms.  We adopt:

\begin{conv}\label{conv:main1}
  If $u$ and $w$ are words in the alphabet $\Sigma$, then $w_u$ will
  denote the number of occurrences of $u$ as a subword of $w$. In
  particular, if $a\in\Sigma$ is a letter, then $w_a$ is the number of
  occurrences of the letter $a$ in $w$ and if  $x,y \in \Sigma$ with $y \ne x^{-1}$
  then $w_{xy}$ is the number of occurrences of $xy$ in $w$.
\end{conv}

\begin{defn}[Whitehead automorphisms]\label{defn:moves}
  A \emph{Whitehead automorphism} of $F_k$ is an automorphism $\tau$ of
  $F_k$ of one of the following two types:

  (1) There is a permutation $t$ of $\Sigma$ such that
  $\tau|_{\Sigma}=t$. In this case $\tau$ is called a \emph{relabeling
    automorphism} or a \emph{Whitehead automorphism of the first
    kind}.

  (2) There is an element $a\in \Sigma$, called the \emph{multiplier}, such
  that for any $x\in \Sigma$
\[
\tau(x)\in \{x, xa, a^{-1}x, a^{-1}xa\}.
\]

In this case we say that $\tau$ is a \emph{Whitehead automorphism of
  the second kind}. (Note that since $\tau$ is an automorphism of $F_k$,
we always have $\tau(a)=a$ in this case). To every such $\tau$ we
associate a pair $(A,a)$ where $a$ is as above and $A$ consists of all
those elements of $\Sigma$, including $a$ but excluding $a^{-1}$, such
that $\tau(x)\in\{xa, a^{-1}xa\}$.  We  say that $(A,a)$ is the
\emph{characteristic pair} of $\tau$.
\end{defn}

Note that for any $a\in \Sigma$ the inner automorphism $ad(a)$ is a
Whitehead automorphism of the second kind. Observe also that the set $SM$ of strictly minimal words is closed under applying relabeling Whitehead
automorphisms, cyclic permutations and taking inverses.

The following is an immediate corollary of Proposition~\ref{wh}.

\begin{prop}\label{easy}
  Let $w$ be a cyclically reduced word of length $n>0$ such that $w\in
  SM$. Let $w'$ be a cyclically reduced word of length $n$. 

Then $w'\in
  Aut(F_k)w$ if and only if there is a relabeling Whitehead automorphism
  $\tau$ such that $w'$ is a cyclic permutation of $\tau(w)$.
\end{prop}

\begin{rem}\label{prim}
  It is easy to see that primitive elements of $F_k$ are never strictly
  minimal.

  If $u \in F_k$ is primitive and $|u|>1$ then $u$ is not minimal and
  hence not strictly minimal. Suppose now that $|u|=1$, so that $u$ is
  $a_i^{\epsilon}$ (where $\epsilon\in \{ 1, -1\}$). Pick an index
  $j\ne i$, $1\le i\le j$. Consider the Whitehead automorphism $\tau$
  of the second kind which sends $a_j$ to $a_ja_i$ and fixes all $a_t$
  for $t\ne j$. Then $\tau(u)=u$, and hence $u$ is not strictly
  minimal. 
\end{rem}

\begin{defn}[Weighted Whitehead graph]

  Let $w$ be a nontrivial cyclically reduced word in $\Sigma^*$.  Let
  $c$ be the first letter of $w$. Thus the word $wc$ is freely
  reduced. (We  use the word $wc$ so that we need only consider
  linear words as opposed to cyclic words.)

  The \emph{weighted Whitehead graph $\Gamma_w$ of $w$} is defined as
  follows.  The vertex set of $\Gamma_w$ is $\Sigma$. For every
  $x,y\in \Sigma$ such that $x\ne y^{-1}$ there is an undirected edge
  in $\Gamma_w$ from $x^{-1}$ to $y$ labeled by the sum $\hat
  w_{xy}:=(wc)_{xy}+(wc)_{y^{-1}x^{-1}}$. where $(wc)_{xy}$ is the number of
  occurrences of $xy$ in $wc$ and $(wc)_{y^{-1}x^{-1}}$ is the number of
  occurrences of $y^{-1}x^{-1}$ in $wc$.

\end{defn}

One can think of $\hat w_{xy}$ as the number of occurrences of $xy$
and $y^{-1}x^{-1}$ in the ``cyclic'' word defined by $w$.  There are
$k(2k-1)$ undirected edges in $\Gamma_w$. Edges may
have label zero, but there are no edges from $a$ to $a$ for
$a\in \Sigma$.  It is easy to see that for any cyclic permutation $v$
of $w$ or of $w^{-1}$ we have $\Gamma_{w}=\Gamma_v$.

\begin{conv}
  Let $w$ be a fixed nontrivial cyclically reduced word.  For two
  subsets $X,Y\subseteq \Sigma$ we denote by $X.Y$ the sum of all
  edge-labels in the weighted Whitehead graph $\Gamma_w$ of $w$ of
  edges from elements of $X$ to elements of $Y$. Thus for $x\in
  \Sigma$ the number $x.\Sigma$ is equal to $w_{x}+w_{x^{-1}}$, the
  total number of occurrences of $x^{\pm 1}$ in $w$.
\end{conv}

The next lemma, which is Proposition~4.16 of Ch.~I in \cite{LS}, gives
an explicit formula for the difference of the lengths of $w$ and
$\tau(w)$, where $\tau$ is a Whitehead automorphism.

\begin{lem}\label{lem:LS}
  Let $w$ be a nontrivial cyclically reduced word and let $\tau$ be a
  Whitehead automorphism of the second kind with the characteristic
  pair $(A,a)$. Let $A'=\Sigma-A$. Then
 \[  ||\tau(w)||-||w||=A.A'-a.\Sigma. \]
\end{lem}

Proposition~\ref{easy} guarantees fast performance of 
Whitehead's algorithm on strictly minimal words.
It turns out that a cyclically reduced word $w$ is strictly minimal if the distribution of the numbers on the edges
of the weighted Whitehead graph of $w$, divided by $|w|$, is close to
the uniform distribution as are the frequencies with which individual
letters occur in $w$.

\begin{lem}[Strict Minimality Criterion]\label{LS}

  Let $0<\epsilon< \frac{2k-3}{k(2k - 1)(4k - 3)}$. Suppose $w$ is a
  cyclically reduced word of length $n$ such that:

  a) For every letter $x\in \Sigma$ we have $\frac{w_x}{n}\in
  (\frac{1}{2k}-\frac{\epsilon}{2},\frac{1}{2k}+\frac{\epsilon}{2})$.

  b) For every edge in the weighted Whitehead graph of $w$ the label
  of this edge, divided by $n$, belongs to
  $(\frac{1}{k(2k-1)}-\epsilon, \frac{1}{k(2k-1)}+\epsilon)$.

  Then for any non-inner Whitehead automorphism $\tau$ of
  $F(a_1,\dots, a_k)$ of second kind we have $||\tau(w)||>||w||=|w|$,
  so that $w\in SM$.
\end{lem}

\begin{proof}

  Let $(A,a)$ be the characteristic pair of $\tau$ and let
  $A'=\Sigma-A$.  Since $\tau$ is assumed to be non-inner, we have
  both $|A| \ge 2$, and $|A'| \ge 2$. Hence $|A|\ |A'| \ge 2(2k-2)$
  and there are at least $2(2k-2)$ edges between $A$ and $A'$ in the
  weighted Whitehead graph of $w$. Recall that $a. \Sigma$ is the
  total number of occurrences of $a^{\pm 1}$ in $w$.

  By Lemma~\ref{lem:LS}, $||\tau(w)||-||w||=A.A'-a.\Sigma$.  By
  assumption on $w$ we have $a.\Sigma\le n(\frac{1}{k}+\epsilon)$ and
\[
||\tau(w)||-||w||=A.A'-a.\Sigma \ge
2n(2k-2)(\frac{1}{k(2k-1)}-\epsilon)- n(\frac{1}{k}+\epsilon)>0,
\]
where the last inequality holds by the choice of $\epsilon$.
\end{proof}

We will see later that the Strict Minimality Criterion holds for an
exponentially generic set of cyclically reduced words.

\section{A little probability theory}\label{Sect:probab}

Fortunately, probability theory provides us with a good way of
estimating the relative frequencies with which particular one- and
two-letter words occur as subwords in freely reduced words of length
$n$ in a free group $F_k$. This tool is called ``Large Deviation
Theory''.  Since we are only interested in applications of Large
Deviation Theory, we refer the reader to Chapter 3 of the excellent and
comprehensive book of Dembo and Zeitouni~\cite{DZ}
on the subject and  give only a brief overview of how
this theory works. The  statements most relevant to our discussion are
Theorem~3.1.2, Theorem~3.1.6 and Theorem~3.1.13 of \cite{DZ}.

\begin{conv}\label{P}
  Let $\Sigma$ be as in Convention~\ref{conv:main}.  Suppose
  $\Pi=(\Pi_{ij})_{i,j\in \Sigma}$ is the transition matrix of a
  Markov process with a finite set of states $\Sigma$. Suppose $\Pi$
  is \emph{irreducible}, that is,  for every position $(i,j)$ there is
  $m>0$ such that $(\Pi^m)_{i,j}>0$. Assume also that $\Pi$ is
  \emph{aperiodic}, that is,  for each $i\in \Sigma$ the $gcd$ of all
  $m>0$ such that $(\Pi^m)_{i,i}>0$ is equal to $1$.
Suppose also that the Markov
  process starts with some probability distribution on $\Sigma$.  Let
  $f:\Sigma\to \mathbb R$ be a fixed function.  Let $Y_1,\dots,
  Y_n,\dots$ be a Markov chain for this process.  We are interested in
  estimating the probability that $\frac{1}{n}\sum_{i=1}^n f(Y_i)$
  belongs to a particular interval $J\subseteq \mathbb R$, or, more
  generally, to a particular Borel subset of $\mathbb R$.  This
  probability defines what is referred to as an \emph{empirical
    measure} on $\mathbb R$. A similarly defined \emph{pair empirical measure} counts $\frac{1}{n}\sum_{i=1}^n g(Y_i,Y_{i+1})$, where $g: \Sigma\times \Sigma\to \mathbb R$ is some function (in the summation one takes $Y_{n+1}=Y_1$).
\end{conv}

\begin{exmp}\label{ex:free}
  In a typical application to free groups, a freely reduced word
  $w=Y_1\dots Y_n$ in a free group $F(a_1,\dots, a_k)$, $k>1$, can be
  viewed as such a Markov chain for a Markov process with the set of
  states $\Sigma=\{a_1\dots, a_k, a_1^{-1},\dots, a_k^{-1}\}$, and
  with transition probabilities $\Pi_{x,y}=P(x|y)=\frac{1}{2k-1}$ if
  $y\ne x^{-1}$ and $\Pi_{x,y}=P(x|y)=0$ if $y=x^{-1}$, where $x,y\in
  \Sigma$. The initial distribution on $\Sigma$ is uniform, so that
  for any $x\in \Sigma$ the probability for a Markov chain to start at
  $x$ is $\frac{1}{2k}$.  The sample space for the Markov process of
  length $n$ consists of \emph{all} words of length $n$ in $\Sigma$.
  However, a word which is not freely reduced will occur as a
  trajectory with zero probability because of the definition of
  $\Pi_{x,y}$. It is easy to see that this Markov process induces
  precisely the uniform distribution on the set of all freely reduced
  words of length $n$ and the probability assigned to a freely reduced
  word of length $n\ge 1$ is $\frac{1}{2k(2k-1)^{(n-1)}}$.

  If we want to count the number $w_a$ of occurrences of $a\in \Sigma$
  in such a freely reduced word, we should take $f$ to be the
  characteristic function of $a$, that is $f(a)=1$ and $f(y)=0$ for
  all $y\ne a$, $y\in \Sigma$.  Then $\frac{1}{n}\sum_{i=1}^n f(Y_i)$
  is precisely $\frac{w_a}{n}$. Similarly, if $g(a,b)=1$ and $g(x,y)=0$ for $(x,y)\ne (a,b)$ then the pair empirical measure essentially counts $\frac{w_{ab}}{n}$.
\end{exmp}

Going back to the general case, Large Deviation Theory guarantees the
existence of a \emph{rate function} $I(x)\ge 0$ (with some additional
good convexity properties) such that for any closed subset $C$ of
$\mathbb R$:

\begin{equation}\label{LDP}
\limsup_{n\to\infty} \frac{1}{n}\log
P(\frac{1}{n}\sum_{i=1}^n f(Y_i)\in C)\le -\inf_{x\in C} I(x).
\end{equation}
Therefore, if $\inf_{x\in C} I(x)=s>0$ then for all but finitely many
$n$ we have
\[
P(\frac{1}{n}\sum_{i=1}^n f(Y_i)\in C)\le exp(-sn/2)
\]
and thus the above probability converges to zero exponentially fast
when $n$ tends to $\infty$.

Similarly, for any open subset $U\subseteq \mathbb R$ we have
\[
\liminf_{n\to\infty} \frac{1}{n}\log P(\frac{1}{n}\sum_{i=1}^n
f(Y_i)\in U)\ge -\inf_{x\in U} I(x),
\]
so that for $s'=\inf_{x\in U} I(x)\ge 0$ we have
\[
P(\frac{1}{n}\sum_{i=1}^n f(Y_i)\in U)\ge exp(-2s'n)\tag\text{\dag}
\]
for all sufficiently large $n$.

Large Deviation Theory also provides an explicit formula for computing
the rate function $I(x)$ above and assures that in reasonably good
cases, like Example~\ref{ex:free} above, the function $I(x)$ is a
strictly convex non-negative function achieving its unique minimum at
a point $x_0$ corresponding to the expected value of $f$ (or the
``equilibrium'').  For instance, in the case of the Markov process for
$F_k$ considered in Example~\ref{ex:free}, the symmetry
considerations imply that $x_0$ is the expected value of the number of
occurrences of $a=\in \Sigma=\{a_1,\dots, a_k, a_1^{-1}, \dots,
a_k^{-1}\}$, divided by $n$, in a freely reduced word $w$ of length
$n$ in $F_k$, that is $x_0=\frac{1}{2k}$. Then
$I(x_0)=0$ and Large Deviation Theory (namely Theorem~3.1.2,
Theorem~3.1.6 of \cite{DZ}) implies that for any $\epsilon>0$ we have
\[\inf \{I(x) | x\in [0, \frac{1}{2k}-\epsilon]\cup
[\frac{1}{2k}+\epsilon, 1]\}=s_{\epsilon}>0.\]

The above computation means that for any fixed $\epsilon>0$ we
have

\begin{gather*}
P(\frac{w_a}{n}\in [0, \frac{1}{2k}-\epsilon]\cup
[\frac{1}{2k}+\epsilon, 1] | w\in F_k \text { with
} |w|=n)\underset{n\to\infty}{\le}\\
\le exp(-s_{\epsilon} n/2),
\end{gather*}
that is, the above probability tends to zero exponentially fast when
$n$ tends to infinity.

Accordingly,
\[
P(\frac{w_a}{n}\in (\frac{1}{2k}-\epsilon,\frac{1}{2k}+\epsilon)| w\in
F_k \text { with } |w|=n) \to_{n\to\infty} 1
\]
and the convergence is exponentially fast.

We  present a formula for computing $I(x)$ for reference purposes.
Let $\Pi, \Sigma, f$ be as in Convention~\ref{P}.  Then formula
\eqref{LDP} holds with

\[
I(x)=\sup_{\theta\in R} \theta x - \log \rho(\Pi_{\theta}).
\]
Here $\Pi_\theta$ is a $\Sigma\times \Sigma$-matrix, where the entry
in the position $(i,j)$ is $\Pi_{ij} exp( \theta f(j))$ and where
$\rho(\Pi_{\theta})$ is the Perron-Frobenius eigenvalue of
$\Pi_\theta$. The convexity of $I(x)$ follows from the fact that in
the above formula $I(x)$ is obtained via a Legendre-Fenchel transform
(also known as ``convex conjugation'') of a smooth function.
A different explicit formula for $I(x)$ is given in
Theorem~3.1.6 of~\cite{DZ}

Dembo and Zeitouni~(see Theorem~3.1.13 of \cite{DZ}) also provide an
analogue of \eqref{LDP} for the pair empirical measure corresponding
to a finite state Markov process, which, in the context of
Example~\ref{ex:free} allows one to estimate the expected relative
frequencies with which a fixed two-letter word occurs as a subword of
a freely reduced word.

Recall that $\gamma(n,F_k)=2k (2k-1)^{n-1}$ is the number of all freely reduced words of length $n$ in $F_k$. When applied to the Markov process corresponding to freely reduced
words in a free group $F_k$, as in Example~\ref{ex:free} above,
Theorem~3.1.2, Theorem~3.1.6 and Theorem~3.1.13 of \cite{DZ} imply the
following:

\begin{prop}\label{markov}
  Let $F_k=F(a_1,\dots,a_k)$ be a free group of rank $k>1$.

  Then:

\begin{enumerate}
\item For any $\epsilon>0$ and for any $a\in \Sigma$ we have
\[
\lim_{n\to\infty} \frac{\#\{w\in F_k \vert \ |w|=n \text{ and }
  \frac{w_a}{n}\in
  (\frac{1}{2k}-\epsilon,\frac{1}{2k}+\epsilon)\}}{\gamma(n,F_k)}=1,
\]
and the convergence is exponentially fast.
\item For any $a,b\in \Sigma$ such that $b\ne a^{-1}$ and for any
  $\epsilon>0$ we have
\[
\lim_{n\to\infty} \frac{\#\{w\in F_k \vert \ |w|=n \text{ and }
  \frac{w_{ab}}{n}\in
  (\frac{1}{2k(2k-1)}-\epsilon,\frac{1}{2k(2k-1)}+\epsilon)\}}{\gamma(n,F_k)}=1,
\]
and the convergence is exponentially fast.
\end{enumerate}
\end{prop}

It is worth noting, as  pointed out to us by Steve Lalley, that one
can also obtain the conclusion of Proposition~\ref{markov} without
using Large Deviation Theory and relying instead on
generating functions methods but  such an approach would be
longer  and  require considerably more computation.

\section{Whitehead graphs of generic words}\label{Sect:genwords}

The following two preliminary statements   are straighforward and we omit the proofs.

\begin{lem}\label{crit}
  The following hold in $F_k$:
\begin{enumerate}
\item For every $n>0$ we have $\gamma(n,C)\le \gamma(n,F_k)\le
  2k\gamma(n,C)$ and $\rho(n,C)\le \rho(n,F_k)\le 2k\rho(n,C)$.
  Moreover,
\[
\gamma(n,F_k)=2k(2k-1)^{n-1} \text{ and }
\rho(n,F_k)=1+\frac{k}{k-1}((2k-1)^{n}-1).\]

\item A set $D\subseteq F_k$ is exponentially $F_k$-negligible if and only
  if $\frac{\gamma(n,D)}{(2k-1)^n}\to 0$ exponentially fast when
  $n\to\infty$.

\item A set $D\subseteq C$ is exponentially $C$-negligible if and only
  if $\frac{\gamma(n,D)}{(2k-1)^n}\to 0$ exponentially fast when
  $n\to\infty$.

\item A subset $D\subseteq F_k$ is exponentially $F_k$-generic if and only
  if
$\frac{\gamma(n,D)}{\gamma(n,F_k)}\to 1$ exponentially fast when
  $n\to\infty$.

\item A subset $D\subseteq C$ is exponentially $C$-generic if and only
  if
$\frac{\gamma(n,D)}{\gamma(n,C)}\to 1$ exponentially fast when
  $n\to\infty$.

\end{enumerate}
\end{lem}

\begin{prop}\label{count}
  Let $A\subseteq C$. Let $A'$ be the set of all freely reduced words
  in $F_k$ whose cyclically reduced form belongs to $A$. Then:
\begin{enumerate}
\item If $A$ is exponentially $C$-negligible then $A'$ is
  exponentially $F_k$-negligible.
\item If $A$ is exponentially $C$-generic then $A'$ is exponentially
  $F_k$-generic.
\end{enumerate}
\end{prop}

The above proposition shows that the notions of being
exponentially $F_k$-generic and exponentially $C$-generic (same for
negligible) essentially coincide.

The results of Large Deviation Theory stated in
Section~\ref{Sect:probab} now allow us to describe the weighted Whitehead
graph of a ``random'' cyclically reduced word of length $n$ of $F_k$.

\begin{prop}\label{use}
  Let $\epsilon>0$ be an arbitrary number.  Let $Q(n,\epsilon)$ be the
  number of all cyclically reduced words $w$ of length $n$ such that
  for every edge of the weighted Whitehead graph of $w$ the label of
  this edge, divided by $n$, belongs to the interval
  $(\frac{1}{k(2k-1)}-\epsilon, \frac{1}{k(2k-1)}+\epsilon)$.
  Similarly, for $a\in \Sigma$ let $T(n,a, \epsilon)$ be the number of
  all cyclically reduced words $w$ of length $n$ such that
  $\frac{w_a}{n}\in (\frac{1}{2k}-\frac{\epsilon}{2},\frac{1}{2k}+\frac{\epsilon}{2})$.

  Then:

\begin{enumerate}
\item We have
\[
\lim_{n\to\infty}\frac{Q(n,\epsilon)}{\gamma(n,C)}=1,
\]
and the convergence is exponentially fast.
\item For any $a\in \Sigma$ we have
\[
\lim_{n\to\infty}\frac{T(n,a,\epsilon)}{\gamma(n,C)}=1,
\]
and the convergence is exponentially fast.
\end{enumerate}
\end{prop}

\begin{proof}
  Denote $N_n=\gamma(n,F_k)$ and $C_n=\gamma(n,C)$.  For a two-letter
  word $xy$ in $\Sigma^*$ denote by $E_{xy}(n,\epsilon)$
  (correspondingly by $E'_{xy}(n,\epsilon)$) the number of all
  cyclically reduced (correspondingly freely reduced) words $w$ of
  length $n$ such that

\[
\frac{w_{xy}}{n}\in [0,\frac{1}{k(2k-1)}-\epsilon] \cup
[\frac{1}{k(2k-1)}+\epsilon,1].
\]
Similarly, for $a\in \Sigma$ let $E_{a}(n,\epsilon)$ (correspondingly
$E'_{a}(n,\epsilon)$) denote the number of all cyclically reduced
(correspondingly freely reduced) words $w$ of length $n$ such that:

\[
\frac{w_{a}}{n}\in [0,\frac{1}{2k}-\epsilon] \cup
[\frac{1}{2k}+\epsilon,1].
\]

Fix a letter $a\in \Sigma$ and a two-letter word $xy$ such that $y\ne
x^{-1}$.

By Lemma~\ref{crit} we know that $C_n\le N_n\le 2k C_n$.  Also, since
every cyclically reduced word is freely reduced, we have
$E_a(n,\epsilon)\le E'_a(n,\epsilon)$ and $E_{xy}(n,\epsilon)\le
E'_{xy}(n,\epsilon)$.

Therefore

\[
\frac{E_a(n,\epsilon)}{C_n}\le 2k \frac{E_a(n,\epsilon)}{N_n}\le 2k
\frac{E'_a(n,\epsilon)}{N_n}\to_{n\to\infty} 0
\]
and
\[
\frac{E_{xy}(n,\epsilon)}{C_n}\le 2k \frac{E_{xy}(n,\epsilon)}{N_n}\le
2k \frac{E'_{xy}(n,\epsilon)}{N_n}\to_{n\to\infty} 0
\]
and the convergence in both cases is exponentially fast by
Proposition~\ref{markov}.

Note that the label, which we denote $\hat w_{xy}$, on the edge
$[x^{-1}, y]$ in the weighted Whitehead graph of a cyclically reduced
word $w$ differs at most by one from $w_{xy}+w_{y^{-1}x^{-1}}$ (since
it is possible that $w$ begins with $y$ and ends with $x$ or that $w$
begins with $x^{-1}$ and ends with $x^{-1}$).

Therefore for all sufficiently large $n$ the condition $|\frac{\hat
  w_{xy}}{n} -\frac{1}{k(2k-1)}|<\epsilon$ implies that
$|\frac{w_{xy}+w_{y^{-1}x^{-1}}}{n} -\frac{1}{k(2k-1)}|<\epsilon/2$.
Let $\hat E_{xy}(n,\epsilon)$ denote the number of all cyclically
reduced words of length $n$ such that $|\frac{\hat w_{xy}}{n}
-\frac{1}{k(2k-1)}|\ge \epsilon$. Then

\[
\frac{\hat E_{xy}(n,\epsilon)}{C_n}\le 2k \frac{\hat
  E_{xy}(n,\epsilon)}{N_n}\le 2k
\frac{E'_{xy}(n,\epsilon/8)+E'_{y^{-1}x^{-1}}(n,\epsilon/8)}{N_n}\to_{n\to\infty}
0
\]
where the convergence is exponentially fast by
Proposition~\ref{markov}. This implies the statement of
Proposition~\ref{use}.
\end{proof}

\section{The generic complexity of Whitehead's  algorithm}\label{Sect:main}

\begin{rem}\label{rem:conj_free}
  Before proving the main result, we need to discuss the complexity of
  the conjugacy problem in the free group $F_k$. Given freely reduced
  words $u',v'$, we can find  their  cyclically reduced forms $u$ and $v$ in 
  time linear in  $\max\{ |u'|,  |v'|\}$) 
  by successively cancelling inverse pairs of letters
  from the two ends of each   word.
  If $|u|\ne |v|$ then clearly $u'$ is not conjugate to $v'$ in $F_k$.

  Suppose now that $|u|=|v|=n$. Then $u'$ is conjugate to $v'$ if and
  only if $u$ is a cyclic permutation of $v$. The naive algorithm 
  of comparing all cyclic permutations of $u$ with $v$ 
  takes quadratic time.
  However,  $u$ is a cyclic permutation of $v$ if and only if $u$ is a
  subword of $vv$.  There is a well-known pattern matching algorithm
  in computer science, called the Knuth-Morris-Pratt algorithm, which decides
  if a word $u$ is a subword of a word $z$ in time linear in $|u| +
  |z|$.  See, for example, \cite{Gus} for details. Applied to the
  words $u$, $vv$, this algorithm allows us to decide if  $u$ is a 
  cyclic permutation of $v$ in linear time in $n$. Thus the conjugacy problem in $F_k$ is actually solvable in 
  time linear in terms of the maximum of the lengths of the two input words.
\end{rem}

We can now prove Theorem~\ref{thm:A} as stated in Section~\ref{main}:

\begin{proof}[Proof of Theorem~\ref{thm:A}]
  Choose $0<\epsilon< \frac{2k-3}{k(2k-1)(4k-3)}$.  Let $L(\epsilon)$
  be the set of all cyclically reduced words $w$ in $\Sigma^*$ such
  that:

  a) for every letter $a\in \Sigma$ we have $\frac{w_a}{n}\in
  (\frac{1}{2k}-\frac{\epsilon}{2},\frac{1}{2k}+\frac{\epsilon}{2})$,
  (where $n=|w|$),

  and

  b) for every edge in the weighted Whitehead graph of $w$ the label
  of this edge, divided by $n$, belongs to
  $(\frac{1}{k(2k-1)}-\epsilon, \frac{1}{k(2k-1)}+\epsilon)$.

  By the Strict Minimality Criterion,  Lemma~\ref{LS}, we have
  $L(\epsilon)\subseteq SM$.  Proposition~\ref{use} and
  Lemma~\ref{crit} imply that $L(\epsilon)$ is exponentially
  $C$-generic. Therefore the bigger set $SM$ is also exponentially
  $C$-generic. Hence by Proposition~\ref{count} the set $SM'$ is
  exponentially $F_k$-generic and part (1) of the theorem is
  established.

  For a fixed Whitehead automorphism $\tau$ and a freely reduced word
  $w\in F_k$ one can compute the freely reduced word $\tau(w)$ in 
  time linear in  $|w|$.  Since the set of Whitehead automorphisms is
  a fixed finite set, one can thus decide in  time linear in  $|w|$ if a
  cyclically reduced word $w$ belongs to $SM$. Thus part (2) of the
  theorem holds. Now Proposition~\ref{wh} together with
  Remark~\ref{rem:conj_free} imply part (3), since there are only
  finitely many relabeling Whitehead automorphisms of the first kind.

  In turn part (3) together with Proposition~\ref{wh} implies parts
  (4) and (5).

\end{proof}

\begin{rem}

  As stated in Theorem~\ref{thm:A}, we can indeed decide if a cyclically
  reduced word $w$ is strictly minimal, that is,  $w \in SM$, in
  time linear in  $|w|$ since the number of Whitehead
  automorphisms is fixed and finite. A priori however,  this requires applying
  every Whitehead automorphism of the second kind to $w$ and then
  computing the freely reduced form of the result. This may be 
  undesirable if the rank $k$ of $F_k$ is large  since  the
  number of Whitehead automorphisms of the second kind grows
  exponentially with $k$.

  On the other hand, the  subset $L(\epsilon)$ of $SM$,
  defined as in  the proof of Theorem~\ref{thm:A}
  with  $\epsilon=\frac{2k-3}{2k(2k - 1) (4k - 3)}$,  is still exponentially generic
  according to the Strict Minimality Criterion. 
  The membership problem in $L(\epsilon)$ is solvable much faster.  All we need to do 
  to decide if   $w\in L(\epsilon)$   is to compute the frequencies
  with which the one- and two-letter subwords occur in $w$ and then
  check if they belong to the required intervals. The number of
  the frequencies with which
  one- and two-letter words occur in $w$ only grows quadratically with $k$.
\end{rem}

\section{Stabilizers of generic elements}

The above analysis also allows us to deduce that stabilizers of
generic elements of $F_k$ in $Aut(F_k)$ and in $Out(F_k)$ are very small.

We need to recall the following property of automorphic orbits which
is a direct corollary of Proposition~4.17 in Chapter~I of \cite{LS}.

\begin{prop}\label{LS1}
  Let $w,w'$ be minimal cyclically reduced words with $||w||=||w'||$ and let
  $\alpha\in Aut(F_k)$ be such that $w'=\alpha(w)$. Then there exist
  Whitehead automorphisms $\tau_i$, $i=1,\dots, n$ such that:

\begin{enumerate}
\item We have $\alpha=\tau_n\dots \tau_1$ in $Aut(F_k)$,
\item For each $i=1,\dots, n$ we have $||\tau_i\dots
  \tau_1(w)||=||w||$.
\end{enumerate}
\end{prop}

Recall that $TS$ as the set of all
 $w\in SM$ such that $w$ is not a proper power and such that
for every nontrivial relabeling automorphism $\tau$ of $F_k$ the
elements $w$ and $\tau(w)$ are not conjugate in $F_k$. Also,
$TS'$ is  the set of elements of $F_k$ whose cyclically reduced form
is in $TS$.

It is easy to see that $TS$ is closed under applying relabeling
automorphisms and cyclic permutations.

\begin{lem}\label{stab}
  Let $w\in TS$ be a nontrivial cyclically reduced word. Then:

\begin{enumerate}
\item If $\alpha\in Aut(F_k)$ is such that $\alpha(w)$ is conjugate to
  $w$ then $\alpha$ is an inner automorphism of $F_k$.
\item The stabilizer $Aut(F_k)_w$ of $w$ in $Aut(F_k)$ is the infinite
  cyclic group generated by $ad(w)$.
\item The stabilizer $Out(F_k)_w$ of the conjugacy class of $w$ in
  $Out(F_k)$ is trivial.
\end{enumerate}
\end{lem}

\begin{proof}
  To see that (1) holds, suppose that $w\in TS$ and that $\alpha(w)=w$
  for some $\alpha\in Aut(F_k)$. Recall that $TS\subseteq SM$.
  Proposition~\ref{LS1} and the definition of $SM$ imply that $\alpha$
  is a product $\alpha=\omega\tau$ where $\omega$ is inner and where
  $\tau$ is a relabeling automorphism.  The definition of $TS$ now
  implies that $\tau$ is trivial and hence $\alpha$ is inner, as
  required.

  Parts (2) and (3) follow directly from (1) since the centralizer of
  a nontrivial element $w$  that is not a proper power in $F_k$ is just
  the cyclic group generated by $w$.
\end{proof}

We will show that the set $TS$ is exponentially $C$-generic.

\begin{lem}\label{perm}
  Let $\tau$ be a nontrivial relabeling automorphism of $F_k$. Let
  $B(\tau)$ be the set consisting of all cyclically reduced words $w$
  such that $\tau(w)$ is conjugate to $w$. Then $B(\tau)$ is
  exponentially negligible in $C$.
\end{lem}
\begin{proof}
   We  only  sketch the argument of the proof,  leaving the details to the reader.

  Let $|w|=n>0$ and suppose that $\tau(w)$ is conjugate to $w$, that
  is $\tau(w)$ is a cyclic permutation of $w$. Suppose first that $w$
  is obtained as non-trivial cyclic permutation $\mu$ of the word
  $\tau(w)$. Then $w$ is uniquely determined by its initial segment of
  length $n/2+1$ and by $\mu$. Note that there are at most $n$
  possibilities for $\mu$. Thus the number of such $w$ is bounded
  above by the number $n\gamma(n/2+1,F_k)$ which grows approximately as
  $n(2k-1)^{n/2+1}$ and thus, after dividing by $(2k-1)^n$, tends to
  zero exponentially fast.

  Suppose now that $w=\tau(w)$. Since $\tau$ is induced by a
  nontrivial permutation of $\Sigma$, this implies that $w$ omits at least
  one letter of $\Sigma$. It is easy to see that for each $a\in
  \Sigma$ the set of all cyclically reduced words $w$ with $w_a=0$ is
  exponentially negligible in $C$. This yields the statement of
  Lemma~\ref{perm}.
\end{proof}

\begin{prop}\label{genST}
  The set $TS$ is exponentially generic in $C$.
\end{prop}
\begin{proof}
  Arzhantseva and Ol'shanskii observed~\cite{AO} that the set of
  cyclically reduced words that are proper powers in $F_k$ is
  exponentially $C$-negligible. It is easy to prove this directly by
  an argument similar to the one used in the proof of
  Lemma~\ref{perm}.  Now Lemma~\ref{perm} and the fact that $SM$ is
  exponentially $C$-generic imply that $C-TS$ is contained in a finite
  union of exponentially negligible sets and hence is itself
  exponentially negligible. Therefore $TS$ is exponentially $C$-generic.
\end{proof}

Proposition~\ref{count} implies that the set $TS'$ of all freely
reduced words, whose cyclically reduced form belongs to $TS$, is
exponentially $F_k$-generic.

We summarize the good properties of $TS$ in the following statement
which follows directly from Proposition~\ref{genST}:

\begin{thm}[c.f. Theorem~\ref{thm:B}]\label{genstab}
  We have $TS=TS'\cap C$ and the following hold:
\begin{enumerate}
\item The set $TS$ is exponentially $C$-generic and the set $TS'$ is
  exponentially $F_k$-generic.
\item There is a linear-time algorithm which, given a freely reduced
  word $w$, decides if $w \in TS'$ or if  $w \in TS$).
\item For any nontrivial $w\in TS'$ the stabilizer $Aut(F_k)_w$ of $w$
  in $Aut(F_k)$ is the infinite cyclic group generated by $ad(w)$.
\item For any nontrivial $w\in TS'$ the stabilizer $Out(F_k)_w$ of the
  conjugacy class of $w$ in $Out(F_k)$ is trivial.
\end{enumerate}
\end{thm}

For future use we also need to establish the genericity of the following
set:

\begin{defn}
  Let the set $Z$ consist of all $w\in TS$ such that there is no
  relabeling automorphism $\tau$ such that $\tau(w)$ is a cyclic
  permutation of $w^{-1}$.
\end{defn}

\begin{prop}\label{inv}
  The following hold in $F_k$.
\begin{enumerate}
\item If $w\in Z$ is a nontrivial word then for any $\alpha\in Aut(F_k)$
  we have $\alpha(w)\ne w^{-1}$.
\item The set $Z$ is exponentially $C$-generic.
\end{enumerate}
\end{prop}

\begin{proof}
  Note that by construction the sets $TS$ and $Z$ are closed under
  taking inverses. Let $w\in Z$ be a nontrivial element.

  The definition of $Z$ and Proposition~\ref{LS1} imply that if
  $\alpha(w)=w^{-1}$ for $\alpha\in Aut(F_k)$ then $\alpha$ is a product
  of inner Whitehead automorphisms and hence is inner itself. However
  in a free group a nontrivial element is not conjugate to its
  inverse. This proves (1).

  For a fixed relabeling automorphism $\tau$ let $D(\tau)$ be the set
  of cyclically reduced words $w$ such that $w^{-1}$ is a cyclic
  permutation of $\tau(w)$.

  Thus to see that (2) holds it suffices to show that for each
  nontrivial relabeling automorphism $\tau$ the set $D(\tau)$ is
  exponentially $C$-negligible. The proof is exactly the same as as
  for Lemma~\ref{perm}.  Namely, if $w\in C$, $|w|=n>0$ and $w^{-1}$
  is obtained by a cyclic permutation $\mu$ of $\tau(w)$, then the
  word $w$ is uniquely determined by $\mu$ and by the initial segment
  of $w$ of length $n/2+1$. Since there are $n$ choices for $\mu$, the
  number of such $w$ is bounded by $n\gamma(n/2+1,C)$, which is
  exponentially smaller than $(2k-1)^n$.
\end{proof}

\section{Applications to generic one-relator groups}

We recall the following classical theorem due to Magnus~\cite{Magnus}:

\begin{prop}\label{magnus}
  Let $G=\langle a_1,\dots, a_k | r=1\rangle$ where  $r$ is a
  nontrivial cyclically reduced word in $F_k$.  Let
  $\alpha\in Aut(F_k)$. Then $\alpha$ factors through to an automorphism
  of $G$ if and only if $\alpha(r)$ is conjugate to either $r$ or
  $r^{-1}$ in $F_k$.
\end{prop}

The following surprising result about ``isomorphism rigidity'' of
generic one-relator groups was  obtained by Kapovich and
Schupp~\cite{KS}.

\begin{prop}\label{KS}
  Let $k\ge 2$ and $F_k=F(a_1,\dots, a_k)$. There exists a exponentially
  $C$-generic set $P_k$ of nontrivial cyclically reduced words with
  the following properties:
\begin{enumerate}
\item There is an exponential time algorithm which, given a cyclically
  reduced word $w$, decides whether or not $w\in P_k$.
\item Let $u\in P_k$. Then $G_u$ is an one-ended torsion-free
  word-hyperbolic group and every automorphism of $G_u$ is induced by
  an automorphism of $F_k$. 
\item Let $u\in P_k$ and let $v$ be a nontrivial cyclically reduced
  word in $F_k$. Then the one-relator groups $G_u$ and $G_v$ are
  isomorphic if and only if there exists $\alpha\in Aut(F_k)$ such that
  $\alpha(u)=v$ or $\alpha(u)=v^{-1}$ in $F_k$.
\end{enumerate}
\end{prop}

We now prove  Theorem~\ref{thm:C} stated in Section~\ref{main}.

\begin{proof}[Proof of Theorem~\ref{thm:C}]
  Let $Q_k=P_k\cap Z$, where $P_k$ is from Proposition~\ref{KS}.  The
  set $Z$ is exponentially $C$-generic by Proposition~\ref{inv} and
  the set $P_k$ is exponentially $C$-generic by Proposition~\ref{KS}.
  Hence $Q_k$ is exponentially $C$-generic as the intersection of two
  exponentially $C$-generic sets and part~(1) of Theorem~\ref{thm:C}
  follows from part~(1) of Proposition~\ref{KS}.

  Suppose $u \in P_k$, as in part (2) of Theorem~\ref{thm:C}. Let $\beta$
  be an automorphism of $G_u$. By Proposition~\ref{KS} $\beta$ is
  induced by an automorphism $\alpha$ of $F_k$.
  Proposition~\ref{magnus} implies that $\alpha(u)$ is conjugate to
  either $u$ or $u^{-1}$ in $F_k$. The latter is impossible by
  Proposition~\ref{inv} since $u\in Z$. Thus $\alpha(u)$ is conjugate
  to $u$. Since $u\in TS$, Lemma~\ref{stab} implies that $\alpha\in
  Inn(F_k)$ and hence $\beta\in Inn(G)$. Thus $Aut(G)=Inn(G)$ and
  $Out(G)=1$. Since $G_u$ is non-elementary torsion-free and
  word-hyperbolic, the center of $G_u$ is trivial and so $G_u$ is
  complete.

  Since $G_u$ is torsion-free one-ended word-hyperbolic and $Out(G_u)$
  is finite,  the results of Paulin~\cite{Pau} show that $G_u$ does not admit
  any essential cyclic splittings. By a theorem of
  Bowditch~\cite{Bow} the boundary of $G_u$ is therefore connected and has no
  local cut-points. Since $G_u$ is a torsion-free one-relator group,
  $G_u$ has cohomological dimension two.  Thus $G_u$ is one-ended
  torsion-free hyperbolic of cohomological dimension two and such that
  $\partial G_u$ is connected and has no local cut-points. A theorem of
  Kapovich-Kleiner~\cite{KK} now implies that $\partial G_u$ is
  homeomorphic to either the Menger curve or the Sierpinski carpet
  and, moreover, if the boundary is the Sierpinski carpet then $G_u$
  must have negative Euler characteristic.

  If $k=2$ then the presentation complex of $G_u$ is topologically
  aspherical~\cite{CCH} (since $G_u$ is a torsion-free one-relator
  group) and can thus be used to compute the Euler characteristic of
  $G_u$. The complex has one $0$-cell, two $1$-cells and one $2$-cell
  so that the Euler characteristic of $G_u$ is $1-2+1=0$. This rules out
  the Sierpinski carpet and hence $\partial G_u$ is homeomorphic to the
  Menger curve in this case. This completes the proof of parts (2) and
  (3) of Theorem~\ref{thm:C}.

  Since $Q_k\subseteq TS$, part (4) of Theorem~\ref{thm:C} follows from
  Proposition~\ref{KS} and Proposition~\ref{LS1}.

  By construction the set $Q_k\subseteq TS\subseteq SM$ and
  $Q_k\subseteq P_k$. Now part~(5) of Theorem~\ref{thm:C} follows from
  Proposition~\ref{KS} and Theorem~\ref{thm:A}.
\end{proof}


\begin{thebibliography}{ABC}



\bibitem{AO} G.~Arzhantseva and A.~Ol'shanskii, \emph{Genericity of
    the class of groups in which subgroups with a lesser number of
    generators are free,} (Russian) Mat. Zametki \textbf{59} (1996),
  no. 4, 489--496


\bibitem{A1} G.~Arzhantseva, \emph{On groups in which subgroups with a
    fixed number of generators are free,}(Russian) Fundam. Prikl. Mat.
  \textbf{3} (1997), no. 3, 675--683.

\bibitem{A2} G.~Arzhantseva, \emph{Generic properties of finitely
    presented groups and Howson's theorem,} Comm. Algebra \textbf{26}
  (1998), 3783--3792.

\bibitem{A3} G.~Arzhantseva, \emph{A property of subgroups of infinite
    index in a free group,} Proc. Amer. Math. Soc. \textbf{128}
  (2000), 3205--3210.

\bibitem{Bahls}
P.~Bahls, \emph{A new class of rigid Coxeter groups.}  Internat. J. Algebra Comput.  \textbf{13}  (2003),  no. 1, 87--94.


\bibitem{BMS} A.~Borovik, A.~G.~Myasnikov and V.~Shpilrain, {\it
    Measuring sets in infinite groups}, Computational and Statistical
  Group Theory (R.Gilman et al, Editors), Contemp. Math., Amer. Math. Soc.
 \textbf{298}  (2002), 21--42.


\bibitem{Bow} B.~Bowditch, \emph{Cut points and canonical splittings
    of hyperbolic groups.}  Acta Math.  \textbf{180} (1998), no. 2,
  145--186


\bibitem{BB} 
R. F. Booth, D. Y. Bormotov and A. V. Borovik, 
\emph{Genetic algorithms and equations in free groups and semigroups}, 
Contemp. Math., Amer. Math. Soc. {\bf 349} (2004)



\bibitem{BMMN}
N.~Brady, J.~McCammond, B.~M\"uhlherr and W.~Neumann, \emph{Rigidity of Coxeter groups and Artin groups.} Proceedings of the Conference on Geometric and Combinatorial Group Theory, Part I (Haifa, 2000).  Geom. Dedicata  \textbf{94}  (2002), 91--109


\bibitem{BV} J.~Burillo and E.~Ventura, \emph{Counting primitive
    elements in free groups.}  Geom. Dedicata \textbf{93} (2002),
  143--162

\bibitem{BKL}
J.~Birman, K.~H.~Ko and S.~J.~Lee, \emph{A new approach to the word
and conjugacy problems in the braid groups.} Adv. Math.
\textbf{139} (1998), no. 2, 322--353


\bibitem{Ch94} C.~Champetier, \emph{Petite simplification dans les
    groupes hyperboliques}, Ann. Fac. Sci. Toulouse Math. (6)
  \textbf{3} (1994), no.~2, 161--221.

\bibitem{Ch95} C.~Champetier, \emph{Propri\'et\'es statistiques des
    groupes de pr\'esentation finie}, Adv. Math. \textbf{116} (1995),
  197--262.

\bibitem{Ch00} C.~Champetier, \emph{The space of finitely generated
    groups,} Topology \textbf{39} (2000), 657--680.

\bibitem{Che96} P.-A.~Cherix and A.~Valette, \emph{On spectra of
    simple random walks on one-relator groups,} With an appendix by
  Paul Jolissaint.  Pacific J. Math.  \textbf{175} (1996), 417--438.

\bibitem{Che98} P.-A.~Cherix and G.~Schaeffer, \emph{An asymptotic
    Freiheitssatz for finitely generated groups,} Enseign. Math. (2)
  \textbf{44} (1998), 9--22.


\bibitem{CCH} I.~Chiswell, D.~Collins and J.~Huebschmann,
  \emph{Aspherical group presentations.}  Math. Z.  \textbf{178}
  (1981), no. 1, 1--36

\bibitem{DZ} A.~Dembo and O.~Zeitouni, \emph{Large Deviation
    Techniques and Applications.}  Second edition. Applications of
  Mathematics, 38.  Springer-Verlag, New York, 1998

\bibitem{FGM}
N.~Franco and J.~Gonz\'{a}lez-Meneses, \emph{Conjugacy problem
for braid groups and Garside groups.} J. Algebra \textbf{266}
(2003), no. 1, 112--132

\bibitem{Gar}
F.~A.~Garside, \emph{The braid group and other groups}, Quart. J.
Math. Oxford \textbf{20} (1969),  no. 78, 235--254

\bibitem{Gh} E. Ghys,
\emph{Groupes Al\'eatoires.} Seminar Bourbaki (March 2003), Asterisque, to appear




\bibitem{Grom} M.~Gromov, \emph{Hyperbolic Groups}, in "Essays in
  Group Theory (G.M.Gersten, editor)", MSRI publ. \textbf{8}, 1987,
  75--263


\bibitem{Grom1} M.~Gromov, \emph{Asymptotic invariants of infinite
    groups.}  Geometric group theory, Vol. 2 (Sussex, 1991), 1--295,
  London Math. Soc. Lecture Note Ser., 182, Cambridge Univ. Press,
  Cambridge, 1993


\bibitem{Grom2} M.~Gromov, \emph{Random walks in random groups},
  Geom. Funct. Analysis \textbf{13} (2003), no. 1, 73--146



\bibitem{Gus} D.~Gusfield, \emph{Algorithms on Strings,Trees and
    Sequences}, Cambridge University Press, Cambridge, 1997.

\bibitem{HMM} R.~Haralick, A.~D.~Myasnikov and A.~G.~Myasnikov,
\emph{Heuristics for Whitehead minimization problem}, preprint



\bibitem{KK} M.~Kapovich and B.~Kleiner, \emph{Hyperbolic groups with
    low-dimensional boundary.}  Ann. Sci. \'{E}cole Norm. Sup. (4)
  \textbf{33} (2000), no. 5, 647--669

\bibitem{KMSS} I.~Kapovich, A.~Myasnikov, P.~Schupp and V.~Shpilrain,
  \emph{Generic-case complexity, Decision problems in group theory and
    Random walks}, J. Algebra \textbf{264} (2003), no. 2, 665--694

\bibitem{KMSS1} I.~Kapovich, A.~Myasnikov, P.~Schupp and V.~Shpilrain,
  \emph{Average-case complexity for the word and membership problems
    in group theory}, Advances in Math., to appear

\bibitem{KS} I.~Kapovich and P.~Schupp, \emph{Genericity, the
    Arzhantseva-Ol'shanskii method and the isomorphism problem for
    one-relator groups}, Math. Ann., to appear

\bibitem{KSn} I.~Kapovich and P.~Schupp, \emph{Delzant's $T$-ivariant,
    one-relator groups and Kolmogorov complexity}, preprint, 2003;\\
http://www.arxiv.org/math.GR/0305353

\bibitem{Khan}
B.~Khan,
\emph{The Structure of Automorphic Conjugacy in the Free Group of Rank Two}, Proceedings of the Special Session on Interactions between Logic, Group Theory and Computer Science, Contemp. Math., Amer. Math. Soc. {\bf 349} (2004) 


\bibitem{Lee}
D.~Lee, \emph{Counting words of minimum length in an automorphic
  orbit}, preprint, 2003

\bibitem{LS} R.~Lyndon and P.~Schupp, \emph{Combinatorial Group
    Theory,} Springer-Verlag, 1977. Reprinted in the ``Classics in
  mathematics'' series, 2000.

\bibitem{Magnus} W.~Magnus, \emph{\"{U}ber diskontinuierliche Gruppen
    mit einer definierenden Relation.} J. rein. angew Math. \textbf{163} (1930), 141--165


\bibitem{Mc} J.~McCool, \emph{Some finitely presented subgroups of the
    automorphism group of a free group.}  J. Algebra \textbf{35}
  (1975), 205--213

\bibitem{MM}
A.~D.~Miasnikov and A.~G.~Myasnikov,
\emph{Whitehead Method and Genetic Algorithms,} Contemp. Math., 
Amer. Math. Soc. {\bf 349} (2004), 89--114. 


\bibitem{MS} A.~G.~Myasnikov and V.~Shpilrain, \emph{Automorphic
    orbits in free groups}, J. Algebra \textbf{269} (2003) , no. 1, pp. 18-27

\bibitem{Mos73}
G.~D.~Mostow, \emph{Strong rigidity of locally symmetric spaces.} Annals of Mathematics Studies, No. 78. Princeton University Press, Princeton, N.J.; University of Tokyo Press, Tokyo, 1973

\bibitem{MW}
 B.~M\"{u}hlherr and R.~Weidmann, \emph{Rigidity of skew-angled Coxeter groups.}  Adv. Geom.  \textbf{2}  (2002),  no. 4, 391--415.


\bibitem{Oliv} Y.~Ollivier, 
\emph{Critical densities for random quotients of hyperbolic groups.}  C. R. Math. Acad. Sci. Paris  \textbf{336}  (2003),  no. 5, 391--394.


\bibitem{Ol92} A.~Yu.~Ol'shanskii, \emph{Almost every group is
    hyperbolic}, Internat. J. Algebra Comput. \textbf{2} (1992),
  1--17.


\bibitem{Pau} F.~Paulin, \emph{Outer automorphisms of hyperbolic
    groups and small actions on $R$-trees.}  Arboreal group theory
  (Berkeley, CA, 1988), 331--343, Math. Sci. Res. Inst. Publ., 19,
  Springer, New York, 1991

\bibitem{PrSp}
S.~Prassidis and B.~Spieler, \emph{Rigidity of Coxeter groups.}  Trans. Amer. Math. Soc.  \textbf{352}  (2000),  no. 6, 2619--2642


\bibitem{Ros}
E.~Rosas, \emph{Rigidity theorems for right angled reflection groups.}  Trans. Amer. Math. Soc.  \textbf{308}  (1988),  no. 2, 837--848

\bibitem{Wh} J.~H.~C.~Whitehead, \emph{On equivalent sets of elements
    in free groups}, Annals of Math. \textbf{37} (1936),
  782--800

\bibitem{Z} A.~Zuk, \emph{On property (T) for discrete groups.}
  Rigidity in dynamics and geometry (Cambridge, 2000), 473--482,
  Springer, Berlin, 2002

\end{thebibliography}
\end{document}